\documentclass[11pt,a4paper]{article}

\date{}

\usepackage[a4paper,margin=1in]{geometry}

\usepackage{amssymb,amsmath}
\usepackage{graphicx}
\usepackage{tabularx}     
\usepackage{booktabs}     
\usepackage{array}        
\usepackage{float}        
\usepackage{url}
\usepackage{comment} 
\usepackage{mathrsfs}
\usepackage{algorithm}
\usepackage{algpseudocode}
\usepackage{hyperref}

\usepackage{changepage}
\usepackage{caption}
\usepackage{enumitem}

\usepackage{lmodern} 
\usepackage{comment}
\usepackage{amsthm}
\newtheorem{theorem}{Theorem}

\newtheorem{definition}{Definition}

\def\argmin{\mathop{\rm argmin}\limits}
\def\argmax{\mathop{\rm argmax}\limits}

\usepackage{xcolor}

\def\X{{\cal X}}
\def\Y{{\cal Y}}
\def\reals{\mathbb{R}}

\def\minimize{\mathop{\rm min}\limits}

\def\st{\mathop{s. t.}}

\usepackage{caption}
\begin{document}


\title{A Review of Bilevel Optimization: Methods, Emerging Applications, and Recent Advancements}

\author{Dhaval Pujara and Ankur Sinha}  

\maketitle

\abstract{
This paper presents a comprehensive review of techniques proposed in the literature for solving bilevel optimization problems encountered in various real-life applications. Bilevel optimization is an appropriate choice for hierarchical decision-making situations, where a decision-maker needs to consider a possible response from stakeholder(s) for each of its actions to achieve his own goals. Mathematically, it leads to a nested optimization structure, in which a primary (leader's) optimization problem contains a secondary (follower's) optimization problem as a constraint. Various forms of bilevel problems, including linear, mixed-integer, single-objective, and multiobjective, are covered. For bilevel problem solving methods, various classical and evolutionary approaches are explained. Along with an overview of various areas of applications, two recent considerations of bilevel approach are introduced. 
The first application involves a bilevel decomposition approach for solving general optimization problems, and the second
application involves Neural Architecture Search (NAS), which is a prime example of a bilevel optimization problem
in the area of machine learning.
}

\providecommand{\keywords}[1]{%
  \noindent\textbf{Keywords: } #1
}

\section{Introduction}
\begin{figure*}[hbt]
\centering
\includegraphics[width=0.85\textwidth]{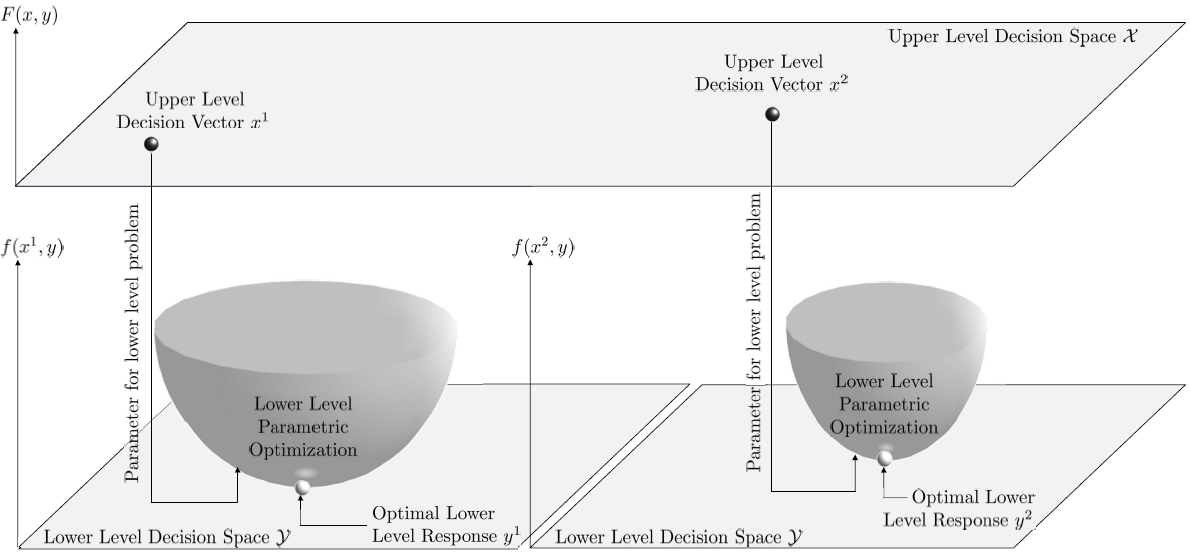}
\caption{A sketch of decision making mechanism in bilevel optimization problem}
\label{fig:bilevel_structure}
\end{figure*}

Bilevel optimization is a special class of optimization problems characterized by a unique structure, where the primary optimization problem contains an additional optimization problem, i.e., a secondary optimization problem, as one of its constraints. In the literature \cite{sinha2017review,dempe02}, the primary and secondary optimization problems are referred to as the upper level and lower level optimization problems, respectively. From a game theory point of view, this setup represents a hierarchical decision-making scenario, where two entities, leader (upper level) and follower (lower level), are associated with each other in a way that the leader needs to think of every possible response from the follower for each of its strategies for achieving its own goals. Thus, the leader's action depends on the reaction from the follower. This mechanism is shown in Figure \ref{fig:bilevel_structure}, representing the inter-linkage between upper and lower levels. The leader has multiple choices or strategies from the upper level decision space $\X$. The decision vector associated with any strategy is represented as an upper level decision vector $x^p$ ($p=1,...,|\X|$). Leader, aware of the set of possible responses from the follower, considers different decisions ($x^1$, $x^2$,...) from $\X$. Follower responds to each decision of leader ($x^1$, $x^2$,...) with an appropriate decision ($y^1$, $y^2$,...) from the lower level decision space $\Y$. For each upper level decision $x^p$, follower performs a lower level parametric optimization using a suitable optimization method, treating $x^p$ as a parameter. This process efficiently explores the lower level decision space $\Y$ and identifies the best (optimal) lower level response $y^p$ to $x^p$. A pair $(x^p,y^p)$ represents a feasible solution to the upper level optimization problem provided that it satisfies all the constraints in the problem. The overall aim is to find the best solution for upper level optimization problem, a pair $(x^*,y^*)$ that minimizes/maximizes the upper level objective function. Thus, in bilevel optimization problem, the leader (upper level) and the follower (lower level) have their own objectives and constraints and both aim to find optimal solutions while functioning in a described decision-making scenario. Here, we discussed a situation where leader usually has complete knowledge of follower's strategies, while the follower only observes the leader's decision and then reacts optimally. Interestingly, there can be various kinds of uncertainties in bilevel problems, for instance, parameter uncertainty, variable uncertainty and decision making uncertainty. Some of the studies in this direction are \cite{ryu2004bilevel, sinha2015solving, burtscheidt2020bilevel, beck2023survey, lu2016finding, lu2018uncertainty,dhyani2026protection, lozano2024bilevel}.

Bilevel optimization has retained the interest of researchers and practitioners as many real-world situations can be formulated as bilevel optimization problems. For instance, pricing models, network design for supply chain management, and competitors operating in the same market are several examples that naturally fit into the bilevel framework as their problem structure involves hierarchical decision-making. Bilevel approach has also been adopted in policy making, where the central authority aims to achieve Social, Technological, Economic, Environmental, or Political (STEEP) enhancement by controlling the actions of consumers through suitable policy norms (strategies). For example, farmers often overuse fertilizers to increase crop yields, which negatively impacts the environment through land and water pollution. In \cite{my-cec15c,whittaker2016spatial,barnhart2017handling}, authors propose a bilevel model based policy that encourages farmers to reduce fertilizer usage and indirectly prevents pollution. According to policy norms, government provides an incentive to farmers if they use fertilizers within a specific range. With this incentive, farmers get the same net profit without overusing fertilizers. Thus, bilevel concept is applied in public sector to make positive impact on environment without compromising the interests of farmers. Apart from the policy formation aimed at regulating STEEP factors, bilevel optimization has been widely applied to homeland security problems such as interdiction of nuclear weapons \cite{brown05, fischetti2019interdiction}, border security \cite{brown09, bucarey2021coordinating, casorran2019study}, defending terror attacks \cite{ordonez2012deployed, shan2013hybrid, baggio2021multilevel}, and protecting critical infrastructure \cite{wein09, an13, furini2019maximum}. In recent times, computer scientists are using bilevel programming for tuning the hyperparameters of various algorithms. Hyperparameters are configuration variables that determine the structural and learning characteristics of an algorithm, e.g., in neural networks, hyperparameters include the learning rate, number of layers, number of neurons per layer, and batch size; in K-means clustering, the number of clusters; in decision trees, the maximum depth and minimum number of samples per leaf; and in evolutionary algorithms, the population size, crossover and mutation probabilities. The values of hyperparameters are set before starting the model training process. Although hyperparameter values are typically set by the user, several studies \cite{bennett2008bilevel, sinha2020gradient, my-gecco14, liu2022bome, petrulionyte2024functional} have developed bilevel based methods to perform the same task. 

From a problem-solving perspective, the hierarchical decision-making structure often leads to non-convex and disconnected feasible regions, which imparts the NP-hard property to bilevel problems and make them difficult to solve mathematically \cite{hansen1992new, vicente1994descent}. A study by Deng \cite{deng1998complexity} provides a proof that no polynomial-time algorithm exists for solving linear bilevel optimization problems. Due to these challenges, traditional mathematical programming based optimization methods fail to solve complex bilevel optimization problems efficiently. 
Apart from that, metaheuristic algorithms (e.g., genetic algorithm, simulated annealing, etc.) are found to be effective in handling some of the inherent difficulties in bilevel problems. Therefore, a combined application of classical methods and metaheuristic algorithms has provided promising results for certain classes of challenging bilevel optimization problems \cite{my-ecj10, sinha2013efficient, angelo13, my-ejor17, my-joh20}. 
Overall, the literature suggests various approaches and methods for solving bilevel problems. Nevertheless, due to the difficulties exhibited by bilevel problems, there is a significant scope for developing new methods to address these problems more efficiently.

\begin{figure*}[hbt]
\centering
\includegraphics[width=0.85\textwidth]{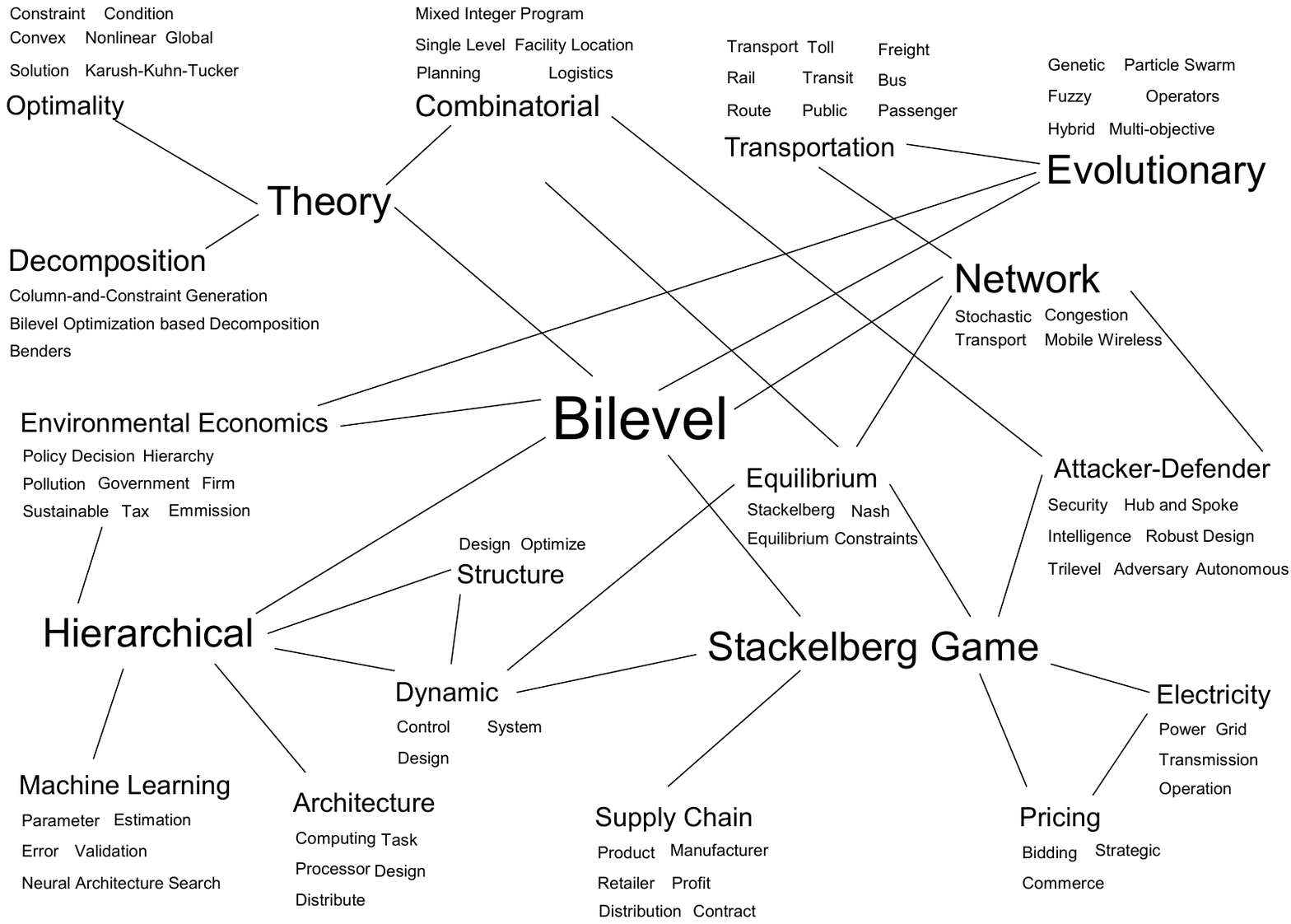}
\caption{A network map depicting connections between research topics addressed using bilevel optimization}
\label{fig:network_map}
\end{figure*}

A network map in Figure \ref{fig:network_map} depicts the applied and theoretical research topics addressed using bilevel optimization and interconnections among them. Each link in map connects a subtopic to a higher level topic represented by a larger font size. 
Along with bilevel optimization methods, this paper provides a review of recent advancements in bilevel optimization research and its application to real-world problems. Accordingly, we start with a basic definition and mathematical presentations of bilevel optimization problem along with the most frequently used terminologies in bilevel optimization. Later, the classical and evolutionary approaches for solving bilevel optimization problems are explained. The method sections are followed by mixed integer bilevel optimization and multiobjective bilevel optimization sections. Next section provides a summary of real-life problems that are addressed using bilevel optimization. In the same section, two recent applications of bilevel optimization approach in (i) solving single-level optimization problems using a bilevel optimization based decomposition method, and (ii) neural architecture search are explained in detail. The paper concludes with a discussion on future research directions in this field.

\section{Bilevel Formulation}
In this section, we formally define the bilevel optimization problem using a mathematical formulation. The notation and setup for bilevel optimization are as follows: the upper level (leader's) optimization problem contains a lower level (follower's) optimization problem nested within it as a constraint. The upper level decision vector is denoted by $x$ ($x=(x_1,...,x_n) \in \reals^n$) and the lower level decision vector is represented by $y$ ($y=(y_1,...,y_m) \in \reals^m$). There are separate sets of objective functions and constraints for the upper level ($F(x,y):\reals^n\times\reals^m \to\reals$ and $G_p(x,y):\reals^n\times\reals^m \to \reals, \,\, p=1,...,P$) and the lower level ($f(x,y):\reals^n\times\reals^m \to\reals$ and $g_q(x,y):\reals^n\times\reals^m \to \reals, \,\, q=1,...,Q$) optimization problems. The lower level problem is a parametric optimization problem to be solved optimally for lower level decision variables with upper level decision variables passed as parameters. For a given $x$, a lower level optimal solution $y$ is considered acceptable at the upper level if $(x,y)$ satisfies all the upper level constraints. The point $(x,y)$ is then a feasible solution to the overall bilevel problem. With this background, a mathematical formulation of the bilevel problem is provided in Definition 1 as follows:

\begin{definition}\label{def:bilevel}
\textit{A bilevel optimization problem with upper and lower level optimization tasks can be formulated as} 
\begin{align}
\mathop{\text{``min''}}_{x,y} \quad & F(x,y) \label{eq:bilevel_start} \\
& \hspace{-16mm} \st &\notag\\
& \hspace{-12mm}y\in \argmin_{y} 
	\lbrace
		f(x,y) : g_q(x,y)\leq 0, \,\, q=1,\dots,Q \rbrace \label{eq:l_level}\\
 & \hspace{-12mm}G_p(x,y)\leq 0, \,\, p=1,\dots,P \label{eq:bilevel_end}
\end{align}
\textit{Above bilevel formulation may also include equality constraints, which are exempted here for brevity. Also, the upper level and lower level decision variables may be integer-valued. However, unless explicitly stated otherwise, we assume these variables to be real and continuous throughout our discussion.}
\end{definition}
The upper level constraints in the presence of both $x$ and $y$ variables are referred to as coupling constraints. It is common in the literature to assume the upper level constraint set as $G_p(x)\leq 0, \,\, p=1,\dots,P$. In the presence of multiple lower level optimal solutions the coupling constraints at the upper level may keep a part of the solutions as feasible while make the other part as infeasible, which leads to various kinds of complexities in defining the problem properly. Therefore, coupling constraints are often avoided in the literature.
The lower level problem (Constraint \ref{eq:l_level}) can be written using a set-valued mapping, which is unknown a priori as follows: 

\begin{definition}\label{def:mapping}
\textit{Let $\Psi:\reals^n \rightrightarrows \reals^m$ be a set-valued mapping,} 
\begin{align}
& \hspace{-11mm} \Psi(x) = \argmin_{y} 
	\lbrace f(x,y) : g_q(x,y)\leq 0, \,\, q=1,\dots,Q \rbrace
\end{align}
\textit{In a nutshell, $\Psi(x)$ represents a mapping that returns the values of lower level variables $y=(y_1,...,y_m)$ for a given $x$. With $\Psi(x)$, the bilevel optimization problem (Definition 1) can be expressed as a general constrained optimization problem as follows:}
\begin{align}
\mathop{\text{``min''}}_{x,y} \quad & F(x,y) \label{eq:mapping_start} \\
& \hspace{-16mm} \st &\notag\\
& \hspace{-12mm} y \in \Psi(x) \\
& \hspace{-12mm} G_p(x,y)\leq 0, \,\, p=1,\dots,P \label{eq:mapping_end}
\end{align}
\end{definition}

In Definitions \ref{def:bilevel} and Definition \ref{def:mapping}, quotes have been used in upper level objective function to reflect the ambiguity that arises in decision-making at upper level when multiple lower level optimal solutions exist for any given upper level decision vector. In this scenario, the decision maker at upper level faces uncertainty regarding which optimal solution will be picked by lower level decision maker. A solution selected from the multiple lower level optimal solutions may or may not benefit the leader, depending on whether the follower's behavior is cooperative or adversarial. The problem becomes fully defined when it is clear which kind of solution will be selected by the lower level in such situations. This aspect is taken into consideration by defining the optimistic and pessimistic positions as follows:

\subsection{Optimistic Position}
In the optimistic position, from the set of multiple lower level optimal solutions, the follower is expected to select a solution that leads to the best objective function value for the upper level or leader's optimization problem. This reflects a certain degree of cooperation between leader and follower. The follower's choice function under the optimistic assumption, $\psi^o(x)$, can be defined as follows:
\begin{align}
& \hspace{-14mm} \psi^o(x) = \argmin_{y} 
	\lbrace F(x,y): y \in \Psi(x) \rbrace \label{eq:optimistic_position}
\end{align}
For simplicity, we have assumed that $ \argmin_{y} \lbrace F(x,y): y \in \Psi(x) \rbrace$ is a singleton and equality sign can be used. Note that if singleton assumption does not hold, one needs to write $$\psi^o(x) \in \argmin_{y} \lbrace F(x,y): y \in \Psi(x) \rbrace$$ and clarify how the choice function picks up the solution.
Accordingly, the formulation of bilevel optimization problem with optimistic choice function (\ref{eq:optimistic_position}) is provided below: 
\begin{align}
\min_{x,y} \quad & F(x,y) \label{eq:bilevel_opt_start} \\
& \hspace{-13mm} \st &\notag\\
& \hspace{-9mm} y = \psi^o(x) \\ \label{eq:bilevel_opt_end}
& \hspace{-9mm} G_p(x,y)\leq 0, \,\, p=1,\dots,P
\end{align}

A bilevel formulation with optimistic position (\ref{eq:bilevel_opt_start})-(\ref{eq:bilevel_opt_end}) is guaranteed to have an optimal solution when it satisfies several mathematical properties mentioned in the theorem below:

\begin{theorem} \label{theo:optimistic}
\textit{If the objective function and constraints of bilevel optimization problem ($F$, $f$, $G_p$, $g_q$) are sufficiently smooth, the constraint region $\Phi$ is non-empty and compact, and the Mangasarian-Fromowitz constraint qualification holds at all points, then the problem is guaranteed to have an optimistic bilevel optimum, provided there exists a feasible solution.}
\end{theorem}

For detailed discussions about the existence of optimistic bilevel optimum and additional results on optimality conditions, readers may refer to \cite{dempe2007new, dempe02, HaPa88, LiMo95, LiMo02, Ou93, dempe2019solution, dempe2019two}.

\subsection{Pessimistic Position}
In the pessimistic position, from the set of multiple lower level optimal solutions, the follower is expected to select a solution that leads to the least favorable outcome for the upper level optimization problem (i.e., leader receives the minimum benefit from the selected lower level optimal solution compared to other lower level optimal solutions). This shows a lack of cooperation between leader and follower. Under pessimistic setting, the choice function of the follower, $\psi^p(x)$, can be defined as follows:
\begin{align}
& \hspace{-14mm} \psi^p(x) = \argmax_{y} 
	\lbrace F(x,y): y \in \Psi(x) \rbrace \label{eq:pessimistic_position}
\end{align}
Once again for simplicity, we assume that $ \argmax_{y} \lbrace F(x,y): y \in \Psi(x) \rbrace$ is a singleton and equality sign can be used. If the singleton assumption does not hold, one needs to write $$\psi^p(x) \in \argmax_{y} \lbrace F(x,y): y \in \Psi(x) \rbrace$$ and clarify how the choice function picks up the solution.
Accordingly, the formulation of bilevel optimization problem with pessimistic choice function (\ref{eq:pessimistic_position}) is provided below: 
\begin{align}
\min_{x,y} \quad & F(x,y) \label{eq:bilevel_pessi_start} \\
& \hspace{-13mm} \st &\notag\\
& \hspace{-9mm} y = \psi^p(x) \\ \label{eq:bilevel_pessi_end}
& \hspace{-9mm} G_p(x,y)\leq 0, \,\, p=1,\dots,P
\end{align}

Between optimistic and pessimistic positions, optimistic position is relatively more tractable as it is possible to reduce the optimistic bilevel formulation to a single-level optimization problem using the variational inequality corresponding to lower level problem, under suitable convexity and regularity assumptions. In the case of pessimistic position, such a straightforward single-level conversion is not possible. As a result, pessimistic bilevel optimization requires explicitly tracking of the lower level optimal solution that yields the worst outcome for the upper level problem, which effectively makes the bilevel problem a three level task. There are certain studies that approximate the pessimistic problems through perturbed optimistic problems, for instance \cite{antoniou2024delta}. Moreover, the pessimistic bilevel problem is guaranteed to have an optimal solution under stronger assumptions provided in the theorem below:

\begin{theorem} \label{theo:pessimistic}
\textit{If the objective function and constraints of bilevel optimization problem ($F$, $f$, $G_p$, $g_q$) are sufficiently smooth, the constraint region $\Phi$ is non-empty and compact, and the set-valued mapping $\Psi(x)$ is lower semi-continuous for all upper level solutions, then the problem is guaranteed to have a pessimistic bilevel optimum.}
\end{theorem}

For details on existence of pessimistic bilevel optimum and additional results on optimality conditions, readers may refer to \cite{dempe02, dempe2014necessary, loridan1996weak, LuMiPi87, wiesemann2013pessimistic, antoniou2025study, dempe2019two, liu2020methods}. For a comprehensive overview of the evolution and advancements in bilevel optimization, readers can refer to review papers \cite{sinha2017review, kleinert2021survey, calvete2020algorithms, smith2020survey, sinha2016evolutionary, dempe2020bilevel_book_chapter, sinha2013bilevel, colson, dempe2003, kalashnikov2015bilevel, ViCa94} and books/edited volumes \cite{bard-book98, dempe02, Co99, ShIsBa97, aksen2013matheuristic, talbi2013metaheuristics, dempe2015bilevel, zemkoho2020bilevel, dempe2020bilevel}. Bilevel optimization is closely related to mathematical program with complementarity constraints (MPCC). For discussion of differences between the two classes of problems, refer to \cite{dempe2012bilevel}.

A summary of the terms and notations commonly used in the bilevel optimization literature is provided in Table \ref{tab:bilevel_summary}.

\begin{table}[htbp] 
\centering
\caption{Key terms and notations used in bilevel optimization literature}
\label{tab:bilevel_summary}
\begin{tabular}{|>{\raggedright\arraybackslash}p{2.8 cm}|
                >{\raggedright\arraybackslash}p{3.0 cm}|
                >{\raggedright\arraybackslash}p{6.6 cm}|}
\hline
\textbf{Terms} & \textbf{Notation(s)} & \textbf{Description} \\
\hline

Decision vectors & 
$x \in X$ \newline $y \in Y$ &
Upper level decision vector ($x$) and decision space ($X$). \newline
Lower level decision vector ($y$) and decision space ($Y$).\\
\hline

Objectives &
$F$ \newline $f$ &
Upper level objective function(s). \newline
Lower level objective function(s). \\
\hline

Constraints &
$G_p, \; p = 1, \ldots, P$ \newline $g_q, \; q = 1, \ldots, Q$ &
Upper level constraint functions. \newline
Lower level constraint functions. \\
\hline

Lower level feasible region &
$\Omega : X \rightrightarrows Y$ &
$\Omega(x) = \{y : g_q(x,y) \leq 0 \; \forall \; q\}$, represents the lower level feasible region for any given $x$.\\
\hline

Constraint region (Relaxed feasible set) &
$\Phi$ &
$\Phi = \{(x,y) : G_p(x,y) \leq 0 \; \forall \; p,\; g_q(x,y) \leq 0 \; \forall \; q \}$, represents the region satisfying both upper and lower level constraints.\\
\hline

Lower level reaction set &
$\Psi : X \rightrightarrows Y$ &
$\Psi(x) = \{y : y \in \operatorname{argmin}_{y \in Y} \{f(x,y) : y \in \Omega(x)\}\}$, shows the lower level optimal solution(s) for a given $x$.\\
\hline

Inducible region (Feasible set) &
$I = \operatorname{gph} \Psi$ &
$I = \{(x,y) : (x,y) \in \Phi, \; y \in \Psi(x)\}$, represents the set of upper level decision vectors and corresponding lower level optimal solution(s) belonging to feasible constraint region. \\
\hline

Choice function &
$\psi : X \rightarrow Y$ &
$\psi(x)$ represents the solution chosen by the follower for any $x$. It becomes important in the case of multiple lower level optimal solutions.\\
\hline

Optimal value function &
$\varphi : X \rightarrow \mathbb{R}$ &
$\varphi(x) = \min\limits_{y} \{f(x,y) : y \in \Omega(x)\}$, represents the minimum lower level function value corresponding to a given $x$.\\
\hline

\end{tabular}
\end{table}

\section{Bilevel Problem Solving Methods}
In the literature, bilevel optimization problems are typically addressed using the classical or evolutionary approaches. Classical approach includes single-level reduction methods, duality methods, descent methods, penalty function methods, and trust-region methods. Broadly, these methods are referred to as mathematical programming based techniques. Single-level reduction methods transform the bilevel problem into single-level optimization problem by replacing the lower level problem with its Karush-Kuhn-Tucker (KKT) conditions, a set-valued mapping, or the lower level optimal value function. Since bilevel problems belong to the class of complex optimization problems, mathematical programming-based techniques are generally applied to bilevel problems that are mathematically well-behaved, typically those with linear, quadratic, or convex objective functions and constraints. Apart from that, strong assumptions such as continuous differentiability and lower semi-continuity are also very common for classical methods. For complex bilevel problems that do not comply with the strong assumptions and favorable mathematical properties (e.g., convexity, continuity, differentiability, etc.), evolutionary approaches are employed, which include nature inspired metaheuristic algorithms such as genetic algorithm, particle swarm optimization, etc. 
In this section, we briefly explain various bilevel problem solving methods covered in the classical and evolutionary approaches.

\subsection{Karush-Kuhn-Tucker Conditions-based Single-Level Reduction} \label{sec:KKT}
In the scenario of lower level problem being convex and sufficiently regular, bilevel optimization problem can be reformulated as a single-level optimization problem by replacing the lower level problem with its Karush-Kuhn-Tucker (KKT) conditions. Accordingly, the problem in Definition \ref{def:bilevel} can be reduced to the formulation given by (\ref{eq:KKT_start})–(\ref{eq:KKT_end}).
\begin{align}
\mathop{\text{min}}_{x,\, y,\, \lambda} \quad & F(x,y) \label{eq:KKT_start} \\
& \hspace{-15mm} \st &\notag\\
& \hspace{-10mm} G_p(x,y)\leq 0, \,\, p=1,\dots,P \\
& \hspace{-10mm} \nabla_y \, L(x,y,\lambda) = 0 \label{eq:lang_const} \\
& \hspace{-10mm} g_q(x,y)\leq 0, \,\, q=1,\dots,Q \\
& \hspace{-10mm} \lambda_q \, g_q(x,y) = 0, \,\, q=1,\dots,Q \label{eq:compl_const} \\
& \hspace{-10mm} \lambda_q \geq 0, \,\, q=1,\dots,Q \label{eq:KKT_end} \\ \notag
& \hspace{-18mm} \textit{where} \\ \notag
& \hspace{-10mm} L(x,y,\lambda) = f(x,y) \, + \sum_{q=1}^{Q} \lambda_q \, g_q(x,y)\\ \notag
\end{align}

The KKT conditions based formulation [(\ref{eq:KKT_start})-(\ref{eq:KKT_end})] sometimes becomes difficult to handle as Lagrangian constraint (\ref{eq:lang_const}) may induce non-convexity, even when the bilevel problem follows convexity and regularity conditions. Also, the complementary conditions (\ref{eq:compl_const}), inherently being combinatorial, make the reduced single-level problem after reformulation a mixed integer program. In the case of linear bilevel problem, the Lagrangian constraint remains linear and complementary conditions are linearized using combinatorial variables, which provides the reduced formulation in Mixed Integer Linear Program (MILP) form. Sometimes, the complementary conditions are also handled using Special Ordered Sets (SOS). To solve the MILP, vertex enumeration \cite{bialas1984two, ChFl92, TuMiVa93} and Branch-and-Bound (B\&B) \cite{FoMc81, BaFa82} approaches are considered in the literature. Though B\&B methods become slow as the number of integer variables increases, these methods are successfully applied for solving single-level reductions of bilevel problems having linear-quadratic \cite{BaMo90} and quadratic-quadratic \cite{AlHoPa92, EdBa91} structures as well. An extended KKT approach has been proposed in \cite{shi2005extended} for linear bilevel problems. For approximate KKT conditions in the context of bilevel programs, readers may refer to \cite{sinha2019using}.

\subsection{Duality-based Single-Level Reduction}
If the lower-level problem is convex and satisfies strong duality, such as through Slater’s condition, the bilevel optimization problem can be equivalently reformulated as a single-level program by replacing the lower-level problem with its dual. This technique is broadly applicable to various classes of bilevel problems that meet the convexity and strong duality criteria.

To illustrate this approach, we will use a general linear bilevel program with continuous variables as an example. However, the underlying principle extends naturally to other convex bilevel formulations where strong duality holds.
\begin{align}
\mathop{\text{min}}_{x, y} \hspace{2mm} & c_{x}^{\top} x+c_{y}^{\top} y \label{eq:linearBilevelStart} \\
\st \notag \\
& A_x x+A_y y \geq a \\
& y \in \argmin_{y} \left\{d^{\top} y: B_x x+B_y y \geq b\right\} \label{eq:linearBilevelEnd}
\end{align}
Next, let us write an alternative single-level reformulation of the bilevel program in (\ref{eq:linearBilevelStart}-\ref{eq:linearBilevelEnd}) using the duality-based approach. 
The dual of the lower level linear program is given as follows:
$$\mathop{\text{max}}_{\lambda} \quad(b-B_x x)^{\top} \lambda \quad \st \quad B_y^{\top} \lambda=d, \lambda \geq 0$$
Now the linear bilevel program in (\ref{eq:linearBilevelStart}-\ref{eq:linearBilevelEnd}), can be reformulated as a single-level problem as follows:
\begin{align}
\mathop{\text{min}}_{x, y, \lambda} \hspace{2mm} & c_{x}^{\top} x+c_{y}^{\top} y \label{eq:linearBilevelDualStart}\\
\st \notag \\  
& A_x x+A_y y \geq a \\ & B_x x+B_y y \geq b \\
& B_y^{\top} \lambda=d, \lambda \geq 0 \\
& d^{\top} y \leq(b-B_x x)^{\top} \lambda \label{eq:linearBilevelDualEnd}
\end{align}
Note that the formulation contains a bilinear constraint, which may not be easy to linearize as the complementary slackness conditions because the bilinear terms do not equate to zero. Refer to \cite{jayaswal2024bilevel} for an overview of using such approaches for practical problem solving.

\subsection{Value Function based Single-Level Reduction}
This approach uses $\varphi$-mapping (optimal value function in Table \ref{tab:bilevel_summary}) to obtain the optimal objective function value for lower level problem $f(x,y)$. This replaces the lower level optimization problem (Constraint \ref{eq:l_level} in Definition \ref{def:bilevel}) with a $\varphi$-mapping constraint and leads to a single-level optimization problem formulation (\ref{eq:phi_start})-(\ref{eq:phi_end}) as follows:
\begin{align}
\mathop{\text{min}}_{x,y} \quad & F(x,y) \label{eq:phi_start} \\
& \hspace{-13mm} \st &\notag\\
& \hspace{-9mm} f(x,y) \leq \varphi(x) \\
& \hspace{-9mm} g_q(x,y)\leq 0, \,\, q=1,\dots,Q \\
& \hspace{-9mm} G_p(x,y)\leq 0, \,\, p=1,\dots,P \label{eq:phi_end}
\end{align}

In practice, the optimal value function $\varphi(x)$ is not known a priori. Therefore, one cannot readily solve the above problem (\ref{eq:phi_start})-(\ref{eq:phi_end}) to arrive at bilevel optimum. Instead, $\varphi$-mapping based algorithms estimate $\varphi(x)$ during the iterations of the algorithm. Some approaches in this direction are \cite{dempe2016solution,lozano2017value,sinha2018bilevel,my-joh20,liu2021value}.

\subsection{Descent Methods}
In bilevel optimization setup, the aim of descent approach is to identify the descent direction that leads to a decrease in the upper level objective function value such that a new point remains feasible. For this feasible point, apart from maintaining the feasibility for upper level problem, it should also be ensuring optimality for the lower level problem. Finding such descent direction can be challenging. To address this, various approaches such as approximating the gradient of upper level objective \cite{KoLa90} and formulating auxiliary problems \cite{SaGa94, ViSaJu94} are suggested in the literature. 

In \cite{SaGa94}, assuming a unique optimal solution, linear independence, convexity, and second-order sufficiency for the lower-level problem, the authors propose solving an auxiliary linear-quadratic bilevel program to determine the steepest descent direction. For the lower level problem,
\begin{align*}
\minimize_{x,y} & \quad f(x,y)\\
\hspace{-5mm} \st \quad  & \\
& \hspace{-8mm} g_i(x,y)\leq 0, i=1,\dots,I
\end{align*}
Let the Lagrangian be represented as follows, where $\mathcal{I}(x) \subseteq \{1,\ldots,I\}$ represents the indices corresponding to the binding constraints at the lower level optimum.
\begin{align*}
L(x,y,\lambda) = f(x,y) + \sum_{i \in \mathcal{I}(x)} \lambda_i g_i(x,y)
\end{align*}
The authors propose solving a linear-quadratic bilevel program to compute the steepest descent direction \( z \in \mathbb{R}^n \) at a point \( x \), given \( y \in \Psi(x) \) and a uniquely determined multiplier vector \( \lambda \) ensured by the assumptions.
\begin{equation*}
\begin{array}{l}
	\minimize_{z, w} \nabla_{x} F\left(x,y\right)^{\top} z + \nabla_{y} F\left(x,y\right)^{\top} w \label{eq:auxULObj}\\
	\st \\
	\begin{array}{ll}
		w \in \argmin_{w} \Big\{ & \left(z^{\top}, w^{\top}\right) \nabla_{(x, y)}^{2} L(x, y,\lambda ) (z, w)\\
		& \st \\
		& \quad \nabla_{y} g_{i}(x, y) w \le -\nabla_{x} g_{i}(x, y) z, \quad i \in \mathcal{I}(x)\\
		& \quad \nabla_{y} f(x, y) w=-\nabla_{x} f(x, y) z+\nabla_{x} L(x, y, \lambda) z \Big\}
	\end{array}\\
 -1 \le z \le 1
\end{array}
\end{equation*}
In the above formulation, the upper level objective function denotes the directional derivative of \( F(x, y) \) along \( (z, w) \), which is minimized to obtain the steepest descent direction. The quadratic program yields \( w \), indicating how the lower-level solution shifts as the upper-level variable \( x \) moves along \( z \).

\subsection{Penalty Function Approach}
Methods from the class of penalty function approach address the bilevel optimization problem by solving a series of unconstrained optimization problems. The unconstrained problem is formed by incorporating a penalty term that measures the violation of the constraints. The penalty term takes the value zero for a feasible solution (i.e., eliminates the penalty term) and takes a positive value (in minimization case) for infeasible solutions (i.e., penalizes the objective function). Penalty function approach was initially implemented in \cite{AiSh81, AiSh84}. Both studies replace the lower level problem with a penalized problem; however, the resultant structure maintains the hierarchy of bilevel optimization, which still remains difficult to solve. Later, \cite{IsAi92} proposed a double penalty method. As the name suggests, this method penalizes both upper and lower level objective functions using the penalty approach. Then, a penalized lower level problem is replaced with corresponding KKT conditions to reduce the bilevel problem to a single-level problem, which is subsequently solved using a penalization technique. There are several studies where lower level optimization problem is directly replaced with corresponding KKT conditions and then penalization approach is applied for solving the single-level problem. In \cite{WhAn93, lv2007penalty}, penalty function approach is considered to solve linear bilevel optimization problems. Former study \cite{WhAn93} converts a bilevel problem to a penalized bilevel problem and solves a series of bilevel problems to achieve the bilevel optimum. Later study \cite{lv2007penalty} performs a single-level reduction using the lower level KKT conditions and updates the upper level objective function by adding complementary slackness conditions with a penalty; then, a reduced single-level problem is solved using series of linear programs.

\subsection{Trust-Region Approach} \label{sec:trust_region}
Trust-region methods perform the local approximation of objective function in the neighborhood of current solution, known as the \textit{trust-region}, where the approximation is assumed to be reliable. These methods are iterative in nature, i.e., build and solve the local approximation of model function in step-by-step manner to reach at the bilevel optimum. 
These methods are found effective for handling non-linearity, non-convexity, or non-regularities in bilevel problems. Trust-region method was first considered in \cite{liu1998trust} to solve non-linear bilevel optimization problem with a lower level problem having convex objective function and linear constraints. 
The study does not contain constraints at the upper level. 
Later, \cite{marcotte2001trust} proposed a more general approach of performing the local approximation of bilevel problem with a linear program at the upper level and a linear variational inequality at the lower level; and then, the solution procedure involves trust-region and line search mechanisms to reach the bilevel optimum over the iterations. Other study \cite{colson2005trust} suggests to approximate the bilevel problem with a linear-quadratic bilevel problem and then solve its reduced single-level formulation as a mixed integer program.

\subsection{Evolutionary Approach} \label{sec:evolu_app}
Evolutionary approach is preferred to deal with complexities such as non-linearity, non-convexity, discontinuity, and non-differentiability in bilevel problems. This approach uses metaheuristic algorithm(s), nature-inspired or intelligence based effective strategies for sampling the solution space, to solve the bilevel problem. Various metaheuristic algorithms used in the literature for solving bilevel optimization problems are genetic algorithm \cite{li07b, mathieu, yin-bilevel, zhu2006hybrid, hejazi2002linear, calvete2011bilevel, li2015genetic, wang2005evolutionary, wang2011new}, particle swarm optimization \cite{li06, jiang2013application, wan2013hybrid}, differential evolution \cite{angelo2014differential, angelo2013differential}, scatter search algorithm \cite{camacho2015heuristic}, etc. Recent comprehensive reviews of metaheuristic approaches for bilevel optimization can be found in \cite{camacho2024metaheuristics,deb2020approximate}. 

Based on the problem solving approach, the bilevel evolutionary methods are classified into three categories: (i) nested methods, (ii) single-level reduction, and (iii) metamodeling. All three types of methods are briefly discussed in this section.

Nested methods solve bilevel problems in nested form (\ref{eq:bilevel_start})-(\ref{eq:bilevel_end}), where lower level problem is solved for each sampled upper level solution. There are two ways considered for solving bilevel problems using nested methods: (i) upper level problem is sampled using metaheuristic algorithm and the corresponding lower level problem is solved using classical method; it is known as hybrid-nested approach and (ii) metaheuristic algorithms are applied at both upper and lower levels to obtain the complete solution. The decision on selecting any approach depends on the complexity of lower level optimization problem. For example, bilevel problems with regular lower level optimization problem are addressed using the hybrid-nested approach \cite{li07b, mathieu, yin-bilevel, zhu2006hybrid}, and the approach of applying metaheuristic algorithm at both levels is considered in the case of bilevel problem with complex lower level optimization problem \cite{li06, my-caor14, camacho2015heuristic, calvete2011bilevel, angelo13, angelo2015study,islam2017enhanced}.

The purpose and mechanism of single-level reduction methods, in the context of evolutionary approach, are similar to what we discussed in Section \ref{sec:KKT}, where bilevel problem is reduced to a single-level problem by replacing the lower-level problem with its KKT conditions, provided that the lower level problem satisfies certain regularity conditions. Most of the time, it is observed that reduced single-level problems continue to sustain various complexities, because of which solving the reduced single-level problem is also not a straightforward task. In such a scenario, evolutionary approach is useful due to better capability of metaheuristic algorithms in handling non-regularities. For example, one of the earliest studies based on evolutionary approach, \cite{hejazi2002linear}, applies the single-level reduction over linear bilevel problem and solves the single-level problem using genetic algorithm, where chromosomes emulate the vertex points. Another study \cite{wang2008genetic} uses a simplex-based genetic algorithm to solve single-level formulation corresponding to linear-quadratic bilevel problem. Other approaches that use evolutionary algorithms with KKT-based reduction and often rely on additional optimization principles are \cite{sinha2019using, li2015genetic, wan2013hybrid, wang2011new, wang2005evolutionary}.

Metamodeling approach is used in optimization when actual function evaluation is very time-consuming or computationally expensive. A meta-model or surrogate model is an approximation of the original model and it is relatively fast to evaluate. Solving a bilevel optimization problem, nested formulation (\ref{eq:bilevel_start})-(\ref{eq:bilevel_end}), leads to a large number of evaluations as we solve lower level problem for each upper level solution, which is a computationally expensive task. If the lower level optimization problem is complex, then this task becomes even more expensive. In this scenario, a combined application of meta-model and metaheuristic algorithm has been observed as an efficient strategy. Reaction set mapping and optimal lower level function value, $\Psi$-mapping and $\varphi$-mapping provided in Table \ref{tab:bilevel_summary}, are often approximated with meta-models. For any given upper level solution vector $x$, $\Psi$-mapping provides the lower level decision vector and $\varphi$-mapping returns the optimal value of lower level optimization problem. In general, none of the mappings are available at the beginning of bilevel problem solving. Therefore, initially, lower level problem is solved for a few upper level solutions, and later the required mapping is approximated ($\hat{\Psi}(x)$ or $\hat{\varphi}(x)$) using the lower level problem elements ($y$ or $f$) and the corresponding upper level solutions. For complex lower level optimization problems, it is hard to approximate the entire mapping, hence, practitioners consider the iterative meta-modeling approach, where required mapping is approximated locally several times over the iterations. After approximating the lower level problem with $\hat{\Psi}(x)$ or $\hat{\varphi}(x)$ mapping, the reduced single-level problem is solved using metaheuristic algorithm. For evolutionary approach where meta models are used, readers can refer to \cite{sinha2014improved, angelo2014differential, sinha2016solving, my-ejor17, islam2017surrogate, sinha2018bilevel,my-joh20}.

\section{Mixed Integer Bilevel Optimization}
The formulation of mixed integer bilevel optimization problem includes the constraints of a general bilevel problem, as given in (\ref{eq:bilevel_start})-(\ref{eq:bilevel_end}) in Definition \ref{def:bilevel} or (\ref{eq:mapping_start})-(\ref{eq:mapping_end}) in Definition \ref{def:mapping}, along with additional constraints that require one or more variables to take integer values only (i.e., any $x_i \in \mathbb{Z}$ and/or $y_j \in \mathbb{Z}$, where $x_i$ is a component of $x$ and $y_j$ is a component of $y$). One of the earlier studies \cite{vicente1996discrete} on mixed integer bilevel optimization examined the properties of mixed integer linear bilevel programs and conditions for the existence of optimal solutions for the problem configuration. In classical optimization literature, branch-and-bound is one of the most commonly used techniques for handling integer variables in MILP. However, branch-and-bound cannot be directly applied to Mixed Integer Bilevel Linear Problems (MIBLP) due to several challenges in fathoming, as discussed in \cite{moore1990mixed}, which also proposed a branch-and-bound technique with strict fathoming conditions and several heuristics to effectively handle problems with a greater number of integer variables. Later, in a subsequent study, \cite{bard1992algorithm} authors developed an implicit enumeration scheme for mixed integer linear bilevel problems under the condition that integer variables can take binary values only. For the continuous-discrete bilevel problems (where $x \in \reals$ and $y \in \mathbb{Z}$), \cite{dempe1996discrete} proposed a cutting plane method that applies the \textit{Chvatal-Gomory} cut. Studies, such as \cite{caramia2015decomposition, fontaine2014benders, saharidis2009resolution}, use benders-decomposition-based techniques to solve bilevel problems with mixed integer upper level variables and continuous linear lower level variables, while \cite{lozano2017value} uses value-function reformulation for mixed integers at both levels. Fischetti et al. (2017) \cite{fischetti2017new} proposed a general-purpose algorithm for MIBLP, where both upper and lower levels have linear constraints and objective functions, and some/all variables take integer values. In a subsequent study \cite{fischetti2018use}, the authors used intersection cuts to solve MIBLP. In \cite{yue2019projection}, a column and constraint generation-based decomposition algorithm is proposed for the single-level formulation corresponding to MIBLP. 
For branch-and-cut implementations that capitalize on useful cuts, readers can refer to \cite{tahernejad2020branch,denegre2009branch}.
Generalized benders decomposition and an extended KKT transformation are discussed in \cite{grimm2019optimal} for a mixed integer nonlinear multilevel model corresponding to the zonal pricing problem in electricity markets. In addition to classical methods, evolutionary approaches have also been used to solve mixed integer bilevel problems in \cite{aksen2013matheuristic, arroyo2009genetic, camacho2014solving, chaabani2015co, handoko2015solving, hecheng2008exponential, legillon2012cobra, miandoabchi2011optimizing}. For the recent developments in discrete bilevel optimization, readers can refer to \cite{kleinert2021computing, kleinert2019global, grimm2019optimal, avraamidou2019multi, fischetti2019interdiction, kleinert2021closing, liu2021enhanced}.

\section{Multiobjective Bilevel Optimization}
In practical scenario, the leader and/or follower may have multiple objectives. This leads to the general formulation of multiobjective bilevel optimization problem below.
\begin{definition}\label{def:multi_bilevel}
\textit{For the upper level objective function $F(x,y): \reals^n \times \reals^m \rightarrow \reals^s$ and lower level objective function $f(x,y): \reals^n \times \reals^m \rightarrow \reals^t$, the multiobjective bilevel problem is formulated in} (\ref{eq:multi_bilevel_start})-(\ref{eq:multi_bilevel_end}) \textit{as follows:}
\begin{align}
\min_{x,y} \quad & F(x,y)= (F_1(x,y),...,F_s(x,y)) \label{eq:multi_bilevel_start} \\
& \hspace{-15mm} \st &\notag\\
& \hspace{-12mm}y\in \argmin_{y} 
	\lbrace
		f(x,y) = (f_1(x,y),...,f_t(x,y)): \notag \\
        & \hspace{6mm} g_q(x,y)\leq 0, \,\, q=1,\dots,Q \rbrace \label{eq:multi_l_level}\\
 & \hspace{-12mm}G_p(x,y)\leq 0, \,\, p=1,\dots,P \label{eq:multi_bilevel_end}
\end{align}
\textit{In above formulation, $G_p: \reals^n \times \reals^m \rightarrow \reals$, $p=1,\dots,P$ represent the upper level constraints, and $g_q: \reals^n \times \reals^m \rightarrow \reals$, $q=1,\dots,Q$ denote the lower level constraints. Additionally, there may be integer restrictions on upper and/or lower level variables.}
\end{definition}

In the literature, multiobjective bilevel optimization has received very scant treatment compared to single objective bilevel optimization, primarily due to computational and decision-making complexities associated with this class of problems. The optimality conditions for multiobjective bilevel programs are discussed in \cite{gadhi2012necessary, ye2011necessary, bank1983non}. In \cite{shi-xia}, an $\epsilon-$constraint technique is applied at both levels of multiobjective bilevel problem, which results into the $\epsilon-$constraint bilevel problem. The $\epsilon-$parameter is passed by the decision maker, and the problem is solved by replacing the lower level problem with its KKT conditions. Eichfelder \cite{eichfelder, eichfelder2} addressed the multiobjective bilevel optimization problems using a classical approach, where the author used a numerical optimization technique and an adaptive exhaustive search method to solve lower level and upper level problems, respectively. The use of exhaustive search method makes this procedure time-consuming and less effective for the large-scale problems. A number of studies involve multiple objectives at the lower level and a single objective at the upper level, and are referred to as semi-vectorial bilevel optimization problems. A common approach to handle such problems is to scalarize the lower level using parameters $\lambda$, and then consider these parameters as a part of upper level decision vector while converting the lower level to a single objective optimization problem \cite{lv2014solution}. This reduces the overall semi-vectorial problem into a single-objective bilevel optimization problem (i.e. one objective at both levels). This idea has also been exploited for bilevel problems with multiple objectives at both levels \cite{gupta2015evolutionary}. Other papers on semi-vectorial bilevel optimization are \cite{bonnel2006optimality, ankhili2009exact, zheng2011solution, calvete2011linear, lv2014solution, ren2016novel}.

With the popularity of evolutionary algorithms in solving multiobjective (single-level) optimization problems during 1990s and 2000s, in late 2000s, researchers started applying the metaheuristic algorithms for solving multiobjective bilevel problems. Yin \cite{yin-bilevel}  formulated the transportation planning and management problem as multiobjective bilevel program, where upper level has multiobjective and lower level contains a single objective, and solved it using a nested genetic algorithm. Later, Halter and Mostaghim \cite{halter-sanaz} proposed a particle swarm optimization based nested strategy to solve the multiobjective bilevel program representing the chemical system. To deal with several complexities such as non-linearity, non-convexity, and non-differentiability in multiobjective bilevel problems, Sinha et al. \cite{sinha2009towards, sinha2011bilevel, sinha2015towards, sinha2015solving} and Deb \& Sinha \cite{my-cec09a, my-mcdm09, my-ecj10} have suggested evolutionary and hybrid approaches to handle multiple objectives at both the upper and the lower levels. An approximation approach of the set-valued mapping has been performed in the context of multiobjective problems as well \cite{sinha2017approximated}. Many of the above studies in the context of multiple objectives at both levels take an optimistic position from the decision making point of view, i.e., among the Pareto-optimal solutions from the lower level, the upper level decision maker freely chooses the one(s) that are most suitable at the upper level. In many real-world scenarios this may not be an actual scenario. For studies that address such concerns with discussions on optimistic/pessimistic Pareto-optimal frontiers and decision making uncertainties, readers may refer to \cite{antoniou2025study,deb2022minimizing,alves2021new,sinha2016evolutionary,sinha2015solving}.

In \cite{my-ecj10}, authors have also developed a suite of test problems to evaluate and compare the performance of various algorithms designed for solving multiobjective bilevel problems. The latest survey on the main approaches for multiobjective bilevel optimization is available in \cite{mejia2023multiobjective}. Readers may also refer to \cite{linnala12, pieume09, pramnik11, ruuska12, zhang12,zou2021multiobjective,chen2024knee,antoniou2025study} for other developments in multiobjective bilevel optimization.

\section{Real-life Applications of Bilevel Optimization}
Problems arising in various domains such as economics, supply chain management, engineering, and management, among others, often exhibit structures that are well-suited to be addressed using the bilevel optimization approaches. In this section, we briefly describe real-life problems addressed using the bilevel optimization and provide a list of relevant reference studies.
\begin{enumerate}[label=\arabic*)]
    \item \text{Toll Setting Problem:} This problem belongs to the class of network problems. In this problem, the authority, who acts as a leader, wants to optimize the toll rates for the network of roads by considering the behavior of network users, the followers. The studies relevant to this problem are \cite{wang2014bilevel, my-cec15a, migdalas95, marcotte2004bilevel, labbe98, kalashnikov2016heuristic, kalashnikov2010comparison, gonzalez2015scatter, yin2002multiobjective, fan2015optimal, constantin95,brotcorne01}. \vspace{2mm}
    
    \item \text{Environmental Economics:} In this class of problems, the authority wants to tax an organization or individuals who consume a particular commodity to generate the revenue. Excessive consumption of this commodity leads to adverse environmental impacts. Therefore, knowing the consumption behavior at various tax rates, authority aims to control the use of commodity by deciding a tax rate that leads to the prevention of environmental pollution without letting consumers lose too much revenue. Variants of this problem are covered in \cite{amouzegar1999determining, bostian2015incorporating, bostian2015valuing,my-cec13,whittaker2016spatial}. \vspace{2mm}

    \item \text{Interdiction Problems:} Interdiction problems are a class of optimization problems under attack or disruption, where an interdictor (attacker) tries to disrupt or degrade the performance of a system, while a defender (operator) tries to optimize system performance despite the disruption. Interdiction is often considered for the nodes or arcs of a network.
    Depending on who is the leader and who is the follower, interdiction problems are often formulated as attacker-defender or defender-attacker problems. These problems may also lead to multi-level optimization (beyond two levels) when network design and protection is taken into account. Some studies in this direction are \cite{wood1993deterministic,israeli2002shortest,church2004identifying,lim2007algorithms,scaparra2008bilevel,akbari2017tri,ramamoorthy2024exact,nigudkar2026solving}. \vspace{2mm}

    \item \text{Facility or Hub Location:} When deciding on the location of a facility or hub, firms may consider the potential responses of its competitors or customers. This scenario represents the Stackelberg game decision-making environment, which can be effectively addressed using bilevel optimization approach. Location problems under the risk of attacks are commonly addressed as interdiction problems as stated above.
    In \cite{kuccukaydin2011competitive}, study considers a scenario where a firm enters a market by locating a new facility, and competitor responds to that by adjusting the attractiveness of its current facilities. Other facility problems addressed using bilevel optimization can be retrieved from \cite{alekseeva2009hybrid, calvete2013efficient, camacho2014solving, uno2008evolutionary, sun08, maldonado2016analyzing, jin2007bi, caramia2015decomposition, church2004identifying, ramamoorthy2018multiple, reisi2019supply, dan2019competitive, dan2020joint, maric2014metaheuristic, bansal2024capacitated}. \vspace{2mm}

    \item \text{Chemical Industry:} For the chemical process, practitioners often want to decide the state variables and quantity of reactants to achieve the optimal output. In this setup, optimizing the output is the upper level problem and the lower level problem appears as an entropy function minimization problem. This problem is considered in \cite{smith82, raghunathan2003mathematical, clark1990bilevel}. \vspace{2mm}

    \item \text{Optimal Design:} The bilevel approach is frequently employed in topology optimization and structural engineering, where the objective is to determine the optimal shape, material distribution, and quantity of material that minimize the overall weight or cost, or equivalently, maximize the structural strength and stiffness at the upper level. The upper-level problem typically encapsulates design decisions subject to physical and performance constraints, such as limits on displacements, stresses, or contact forces, ensuring structural integrity and feasibility. At the lower level, the problem commonly manifests as a potential energy minimization, compliance maximization, worst-case disturbance identification, or as a variational inequality representing the equilibrium state of the structure under applied loads and boundary conditions. This hierarchical formulation elegantly couples design optimization with mechanical equilibrium, allowing the structural behavior to be consistently reflected in design updates. For this class of problems, readers may refer to \cite{kocvara1995solution, kovcvara1997topology, herskovits2000contact, christiansen2001stochastic, bendsoe95, guo2013robust, pan2023first}. \vspace{2mm}

    \item \text{Inverse Optimal Control:} This problem is mainly observed in computer vision, remote sensing, robotics, and related fields. In control theory, one of the major tasks is to obtain the performance index or reward function that fits best on a given dataset. This task is associated with inverse optimal control theory, where one obtains the calculation of the cause based on the given result. Such requirement requires solving a parameter estimation problem with an optimal control problem. This bilevel nature problem is studied in \cite{albrecht2011imitating, johnson2013inverse, mombaur2010human, my-cec16b}. \vspace{2mm}

    \item \text{Principal-Agent Problems:} It is a classical problem observed in economics area, where a principal (leader) subcontracts a job to an agent (follower). In this problem, it is important for principal to take the agent's preferences into account while designing the incentive scheme as the agent is expected to act in his own interests rather than those of the principal. This setup matches with the bilevel decision making mechanism and can be correlated with real-life situation experienced with doctor-patient, employer-employee, politician-voters, corporate board-shareholders, etc. The studies related to this problem are \cite{Bilevel-linear, xu2007supply,garen1994executive, cecchini2013solving}.  \vspace{2mm}

    \item \text{Energy Networks and Market:} After liberalization of the electricity sector and the introduction of energy markets, private power generation companies, market operators, and transmission system operators have become part of the decision-making process in energy sector. The nature of interaction and sequence in which these entities make decisions matches with a Stackelberg-type environment, i.e., bilevel optimization scenario. Therefore, in recent times, bilevel optimization has been extensively employed to address problems associated with energy networks and markets \cite{dempe2015bilevel}. The vulnerability of power systems and the security of power grids under disruptive threats are discussed in \cite{motto,arroyo2010bilevel}. Problems of power generation, transmission, and capacity planning and expansion are studied in \cite{garces2009bilevel,jin2011capacity,jenabi2013bi}. A recent survey on application of bilevel optimization in electricity market is available in \cite{wogrin2020applications}. Apart from electricity market, energy management in the gas market using bilevel optimization is discussed in \cite{grimm2019multilevel, bottger2022cost, schewe2022global}.

    
\end{enumerate}
Next, we delve into the recent applications of bilevel approach in more detail.

\subsection{Bilevel Optimization based Decomposition}
Bilevel Optimization based Decomposition (BOBD) is a recently developed decomposition method, covered in \cite{sinha2024decomposition, sinha2025bilevel}, for effectively solving the complex and large-scale single-level or general optimization problems. As the name suggests, BOBD method uses a decomposition strategy to transform a general optimization problem into an equivalent bilevel optimization problem. This decomposition allows the use of effective bilevel optimization techniques in solving general optimization problems. The motivation and procedure of formulating a general optimization problem as a bilevel optimization problem and solving it using bilevel optimization methods are discussed in this section. A mathematical formulation of general optimization problem is provided below.
\begin{definition}\label{def:singleLevel}
\textit{The single-level optimization problem, with decision variables $x= (x_1, \dots,x_n)$, objective function $F(x)$, and constraints set $G(x)$, can be formulated as follows:}
\begin{align}
\min_{x} \quad & F(x) \label{eq:startSingle}\\
& \hspace{-9mm} s.t. \hspace{4mm} G_r(x)\leq 0, \quad r=1,\dots,R  \label{eq:endSingle}
\end{align}
\end{definition}

The single-level optimization problem (\ref{eq:startSingle})-(\ref{eq:endSingle}) corresponding to various real-world optimization systems often contains a large number of decision variables and constraints, which leads to a large-scale optimization problem. Additionally, the objective function $F(x)$ and constraints $G(x)$ may contain several terms that lead to non-regularities such as non-linearity, non-convexity, non-differentiability, etc. Solving such complex optimization problems, characterized by large $n$ and $R$ along with non-regularities, using a classical approach or evolutionary approach alone does not yield effective solutions generally. That is because, classical methods can handle the large-scale scenario only when problems follow certain regularity conditions related to linearity, convexity, differentiability, etc. On the other hand, metaheuristic algorithms can effectively handle such irregularities, but their performance often deteriorates as problem size increases. Thus, both classical and evolutionary approaches, when employed independently, fail to handle the large-scale and non-regular scenarios. However, most of the time, real-world optimization problems hold both of these complexities simultaneously. Thus, in practical scenarios, both approaches may become ineffective. BOBD method considers a hybrid approach that allows applying both classical and metaheuristic algorithms simultaneously to solve the bilevel optimization problem corresponding to the single-level optimization problem. The working mechanism of BOBD method is discussed next.

Any constrained single-level optimization problem (Definition \ref{def:singleLevel}) has mainly two components: objective function (\ref{eq:startSingle}) and constraints (\ref{eq:endSingle}). Both of these components are basically the mathematical functions of decision variables $x$, i.e., $F(x)$ and $G(x)$. In other words, $F(x)$ and $G(x)$ are mathematical expressions composed of terms containing decision variables (constants are omitted in this discussion). The complexity of an optimization problem is determined in part by the manner in which its decision variables appear in the objective function and constraints. Certain variables may occur only in mathematically regular terms, whereas others may participate in terms responsible for nonlinearity, nonconvexity, nondifferentiability, or other computational difficulties. To understand this, consider the numerical example provided below.
\begin{align*}
\min_{x} \quad & x_1 - 2x_2^3 + 4x_3\\
& \hspace{-9mm} s.t. \quad -x_1 - x_3 \leq -5;\\
& \hspace{0mm} \quad x_2^2 + x_3 \leq 4;\\
& \hspace{0mm} \quad 0 \leq x_1, \, x_2, \, x_3 \leq 10
\end{align*}
In the above numerical example, term $-2x_2^3$ makes the objective function non-convex and term $x_2^2$ causes non-linearity in the second constraint. Imagining the above example without the terms containing $x_2$ decision variable would turn the problem into a linear optimization problem since all the remaining terms with decision variables $x_1$ and $x_3$ are linear. Thus, variable $x_2$ is the source of complexity in the instance. Based on this analysis, each variable ($x_1$, $x_2$, $x_3$) can be tagged as complexity-causing variable $\{x_2\}$ or complexity-soothing variable $\{x_1,x_3\}$. In this way, current exercise, \textit{complexity analysis}, shows the role of variables in the complexity of problem.

The decomposition strategy of BOBD method is built on the complexity analysis procedure discussed above using an example. It performs complexity analysis for each decision variable $x_i$ ($x_i \in x, \, i\in[1,n]$) and classifies $x_i$ into upper level variables category $u$ ($u \subseteq x$) if $x_i$ is complexity-causing variable; otherwise, when $x_i$ is complexity-soothing variable, $x_i$ is classified into lower level variables category $l$ ($l \subseteq x$). Thus, $u \cup l = x$, with $u \cap l = \emptyset$. This classification further allows us to write objective function $F(x)$ and constraints $G_r(x)$ in terms of upper and lower level variables $u$ and $l$, i.e., $F(u,l)$ and $G_r(u,l)$. Later, based on the presence of upper and lower level variables, each constraint $G_r(u,l)$ is separated as the upper level constraint $G_{r^+}(u,l)$ or the lower level constraint $G_{r^-}(u,l)$. Using this variable and constraint classifications, $u/l$, $G_{r^+}(u,l)/G_{r^-}(u,l)$, $F(u,l)$, the single-level optimization problem (\ref{eq:startSingle})-(\ref{eq:endSingle}) can be represented in the form of bilevel optimization problem as follows:
\begin{align}
\min_{u,l} \quad & F(u,l) \label{eq:startDecom}\\
& \hspace{-10mm} \st \hspace{5mm} l \in \argmin_{l} \lbrace F(u,l) : G_{r^-}(u,l) \leq 0 \rbrace \\
& \hspace{0mm} G_{r^+}(u,l) \leq 0\label{eq:endDecom}
\end{align}
In above formulation, the objective functions at both upper and lower levels are same, i.e., $F(u,l)$. Hence, solving this bilevel formulation (\ref{eq:startDecom})-(\ref{eq:endDecom}) would lead to the same optimal solution as the formulation (\ref{eq:startSingle})-(\ref{eq:endSingle}) under certain conditions. If the relevant lower level optima are attained, the bilevel reformulation is equivalent to the original single-level problem provided that every constraint involving at least one lower-level variable is assigned to the lower-level problem, while the constraints retained exclusively at the upper level involve only upper level variables. Since the objective function is identical at both levels, for any fixed upper level decision \(u\), the lower level problem selects an optimal \(l\) over exactly the corresponding feasible choices of the original problem. The upper level optimization over \(u\) therefore recovers the optimal objective value and an optimal solution of the original single-level problem.

For the numerical example considered in ongoing section, an equivalent bilevel formulation can be obtained using a decomposition strategy of BOBD method as follows. Among the three decision variables $x_1$, $x_2$, $x_3$ $\in x$, variable $x_2$, complexity-causing variable, is classified into the upper level variables category (i.e., $x_2=u_{1} \in u$); and variables $x_1$ and $x_3$, complexity-soothing variables, are classified into the lower level variables category (i.e., $x_1,x_3=l_{1},l_{2} \in l$). This classification leads to $x=(x_1,x_2,x_3)= (l_{1},u_{1},l_{2})$. Using this classification, the numerical instance, in the form of single-level optimization problem, can be decomposed into bilevel problem as follows:
\begin{adjustwidth}{1cm}{0cm}
\begin{align*}
\min_{u,l} \quad & l_1 - 2{u_{1}}^3 + 4l_2 \\
& \hspace{-10mm} \textit{s.t.} \\
& l \in \argmin \left\{
\begin{array}{l}
l_1 - 2{u_{1}}^3 + 4l_2 \\
\textit{s.t.} \quad
- l_1 - l_2 \leq -5; \quad {u_{1}}^2 + l_2 \leq 4; \\
\quad \quad \quad 0 \leq l_{1},l_{2} \leq 10
\end{array}
\right\} \\
& 0 \leq u_{1} \leq 10
\end{align*}
\end{adjustwidth}

Next, we focus on the approach considered for solving the derived bilevel optimization problem. Solving an optimization problem simply means determining the values of decision variables such that all the constraints remain satisfied and the maximum/minimum objective function value can be attained when the obtained values of decision variables are inserted into constraints and objective function. The BOBD method classifies the complexity-causing variables into upper level and complexity-soothing variables into lower level; thereafter, an evolutionary approach is used to sample the complexity causing upper-level variables, while a classical method is used to optimize the complexity-soothing lower level variables for each given upper level decision. Therefore, BOBD method follows hybrid approach to solve the bilevel problem. It uses genetic algorithm to sample the values of upper level variables $u$ and solves the lower level optimization problem using mathematical-programming based method to obtain the values of lower level variables $l$ for any given $u$. The stepwise procedure of BOBD method is mentioned in Algorithm \ref{algo:BOBD}.


\begin{center}
\begin{minipage}{0.85\textwidth}   
\begin{algorithm}[H]
\caption{Bilevel Optimization-based Decomposition (BOBD)}\label{algo:BOBD}
\begin{tabularx}{\linewidth}{lX}
\textbf{Input}:  & $F(x)$, $G(x)$- single-level optimization problem\\
\textbf{Output}: & $x^*$- efficient solution of single-level optimization problem\\
\hline \\
\textbf{Step 1}: \,  & Classify each variable $x_i \in x$ into upper level ($u$) or lower level ($l$) categories (i.e., $x=(u,l)$). \vspace{1mm}\\
\textbf{Step 2}: \,  & Formulate a bilevel problem by decomposing the single-level optimization problem with respect to $(u,l)$. \vspace{1mm}\\
\textbf{Step 3}: \,  & Solve the bilevel problem using an evolutionary algorithm at the upper level and a classical algorithm at the lower level (i.e., perform evolutionary sampling of $u$ and use a classical method to solve for $l$ for the given $u$).\vspace{1mm}\\
\textbf{Step 4}: \,  & Return the best obtained solution ($x^*=(u^*,l^*)$) w.r.t. $F(x)$ and $G(x)$.\\
\end{tabularx}
\end{algorithm}
\end{minipage}
\end{center}


To evaluate the performance of the BOBD method, a study by Sinha et al. \cite{sinha2024decomposition} suggested a test suite of 10 test problems (TP1–TP10). Some of these test problems are from a real-world context. Out of 10 test problems, 8 test problems are scalable in terms of variables and constraints and exhibit various complexities such as non-linearity, non-convexity, discontinuity, non-differentiability, etc. The authors have shown the effectiveness of BOBD method by solving each test problem using BOBD method, metaheuristic algorithm (GA), and classical methods (Interior point and Sequential Quadratic Programming) and later comparing the results from each method in small to large-scale scenarios. The BOBD method outperformed all other methods in small, medium, and large-scale scenarios.

\begin{figure*}[h]
\centering
\includegraphics[width=0.85\textwidth]{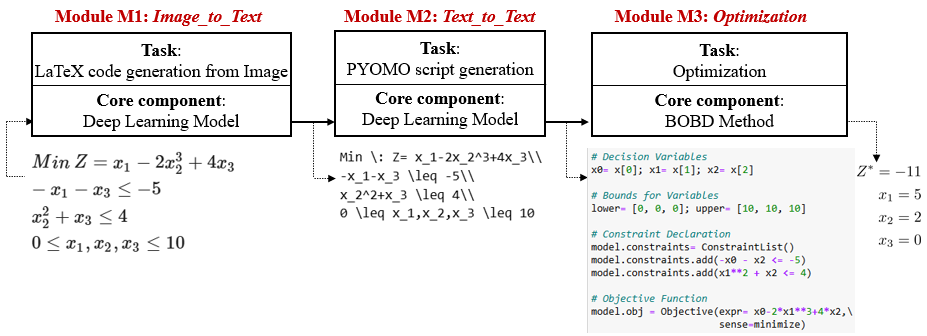}
\caption{AutoOpt framework \cite{sinha_AutoOpt} to automate the optimization problem-solving task}
\label{fig:AutoOpt}
\end{figure*}

In \cite{sinha2024decomposition}, authors manually classify each variable into upper level ($u$) or lower level ($l$) categories (Step 1 in Algorithm \ref{algo:BOBD}). In the next study \cite{sinha2025bilevel}, the authors presented an eigen-value and eigen-vector based variable classification heuristic to automate the variable classification task, which makes BOBD method suitable for large-scale instances by eliminating human intervention. In \cite{sinha_AutoOpt}, authors adopted a machine learning approach to automate the variable classification task, wherein a logistic regression model is built to classify the decision variables into upper and lower level categories. In the same study, BOBD method is part of the \textit{AutoOpt} framework designed to automate the optimization problem-solving task, as shown in Figure \ref{fig:AutoOpt}. \textit{AutoOpt} framework contains three modules in series as follows: M1(\textit{Image\_to\_Text})- contains a deep learning model that considers the image of an optimization formulation and generates the corresponding LaTeX code; M2(\textit{Text\_to\_Text})- contains a deep learning model that extracts the optimization problem from the LaTeX code and generates a corresponding Pyomo script (programming structure for mathematical modeling language); and M3(\textit{Optimization})- contains BOBD method to solve the optimization problem from Pyomo script. The procedure followed by \textit{AutoOpt} framework to solve the example discussed in the current section is shown in Figure \ref{fig:AutoOpt}.

There are several enhancements possible for the BOBD method, for instance, identifying other approaches for variable classification. In BOBD, the complexity-causing variables are classified into upper level and complexity-soothing variables are classified into lower level; hence, lower level optimization problem contains a relatively less complex objective function and constraints. Considering that, several ideas like $\Psi$-mapping, penalty function, or trust-region, can be exploited to reduce the evolutionary sampling at the upper level or the number of calls to the classical algorithm for the lower level.

\subsection{Neural Architecture Search (NAS)}
\begin{figure*}[hbt] 
\centering
\includegraphics[width=0.85\textwidth]{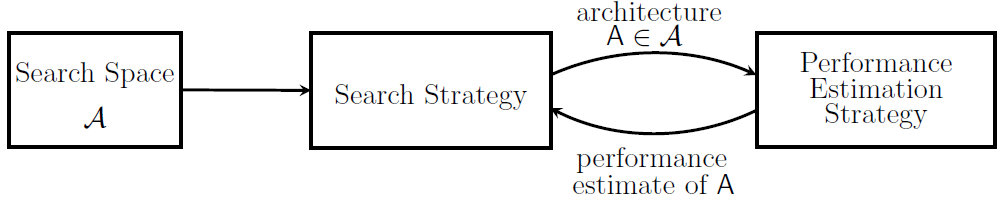}
\caption{An outline of NAS: automated process for identifying the optimal architecture from the complex search space}
\label{fig:NAS_framework}
\end{figure*}

In machine learning (ML), \textit{architecture parameters} refer to the structural parameters of an algorithm that define overall design and configuration of the algorithm. For example, in an \textit{Artificial Neural Network (ANN)}, the number of layers and the number of neurons per layer constitute key architecture parameters. The process of selecting suitable values for these parameters, known as \textit{architecture engineering}, is carried out before training the model on a given dataset. This step is crucial, as the architecture directly influences the model’s learning capacity, generalization ability, and computational efficiency. Traditionally, architecture engineering was performed manually by human experts, making it both time consuming and prone to suboptimal choices. Consequently, there has been growing interest in automating this process through \textit{Neural Architecture Search (NAS)}, which systematically explores possible architectures to identify high-performing designs. It is important to note that architecture engineering and NAS are specialized forms of \textit{hyperparameter optimization}, focusing specifically on the structural aspects of machine learning models.
NAS methods have been reported to outperform manually designed architectures on some tasks such as image classification \cite{zoph2017automl, real2019regularized}, object detection \cite{zoph2017neural}, and semantic segmentation \cite{chen2018searching}. Figure \ref{fig:NAS_framework} illustrates a general mechanism of NAS  \cite{elsken2019neural}. The main aim of NAS is to identify the optimal architecture from a large and complex search space ($\mathcal{A}$) of all possible configurations. Search strategy deals with development of effective technique to explore the search space. Performance estimation strategy dictates how to measure the effectiveness of an architecture ($\textsf{A} \in \mathcal{A}$) obtained from the search strategy.

In this paper, we mainly focus on the use of bilevel optimization approach for neural architecture search (NAS). Accordingly, a mathematical formulation representing the implementation of the general mechanism of NAS (Figure \ref{fig:NAS_framework}) using bilevel optimization approach is provided below:
\begin{align}
\min_{\mathsf{A}} \quad & \mathcal{L}_v(\mathsf{A},W^*) \label{eq:start_NAS}\\
& \hspace{-10mm} s.t. \hspace{5mm} W^* \, \in \, \argmin_{W \in \mathcal{W}} \mathcal{L}_t(\mathsf{A},W) \\
& \hspace{0mm} \textsf{A} \in \mathcal{A} \label{eq:end_NAS}
\end{align}
In (\ref{eq:start_NAS})-(\ref{eq:end_NAS}), the upper level optimization problem aims to minimize the validation loss $\mathcal{L}_v$ w.r.t. the architecture parameters \textsf{A} (upper level variables), and the lower level optimization problem deals with identifying the model parameters $W$ (lower level variables) such that training loss $\mathcal{L}_t$ is minimized for the given \textsf{A}. In this framework, the upper level loop acts as a search strategy whose function is to identify effective architectures by exploring the search space. The lower-level problem represents model training for a fixed architecture, while the upper-level objective evaluates the resulting trained model on the validation set and thereby guides the architecture search. There is an assumption that for any architecture $\textsf{A}^i$ passed by the upper level, there exists at least one optimal model parameters $W^*(\textsf{A}^i)$. This general assumption ensures the well-defined and non-empty feasible region for the upper level problem. Such bilevel approach for NAS enables the joint optimization of architectural choices and model performance. Accordingly, we now discuss several important bilevel approach based NAS methods, which are broadly classified into two categories: (i) sampling-based NAS methods and (ii) bilevel theory-based methods. The core searching mechanisms of methods from both classes are discussed below.
\\
\\
\textbf{Sampling-based NAS Methods}: 
This class of methods uses heuristic or probabilistic sampling strategies to explore the architecture search space. These methods are relatively easy to implement and often serve as a baseline for more sophisticated methods. These methods require a large number of architecture evaluations, which can be computationally expensive. A brief overview of each sampling-based NAS method is provided below.

\begin{enumerate}[label=\arabic*)]
    \item \textbf{Grid Search-based NAS:} This method follows a deterministic approach that involves trying multiple possible combinations of architecture parameters and model parameters. The architecture parameters are defined over a finite and discrete search space.\vspace{2mm}
    
    \item \textbf{Random Search-based NAS:} This method samples the search space using predefined probability distributions. This stochastic approach often yields better results under a fixed computational budget as compared to the grid search method.\vspace{2mm}

    \item \textbf{Bayesian Optimization-based NAS:} This method uses probabilistic surrogate models to identify the promising regions of the architecture search space. It is a more efficient approach compared to grid search and random search approaches. The surrogate models are generally instantiated using Gaussian Processes or Tree-structured Parzen Estimators (TPE). TPE based Bayesian optimization method is widely used for NAS. It constructs a probabilistic model that estimates the likelihood of achieving better performance using a given architecture. The acquisition function, part of the TPE based NAS framework, enables to maintain a balance between exploration and exploitation of search space by strategically selecting the next architecture to evaluate.\vspace{2mm}

    \item \textbf{Evolutionary Computation-based NAS:} Evolutionary Computation-based NAS methods use metaheuristic algorithms to effectively explore complex search spaces without relying on gradient information. Though there are methods that integrate evolutionary exploration with hypergradient-based local search \cite{sinha2023linearCEC}.   
    Population based metaheuristic algorithms, most commonly genetic algorithm or differential evolution, initiate with a population of randomly generated architectures that evolves over iterations as part of the upper level task. At the upper level, crossover and mutation operators are used to update the population or generate new architectures. These architectures are then passed to the lower level, where a model is trained using each architecture and its performance is recorded as the fitness of respective architecture. The latest surveys on the evolutionary computation-based NAS methods are available in \cite{EC_based_NAS_review, kaveh2023application, yao2018taking}.\vspace{2mm}  
  
    \item \textbf{Reinforcement Learning-based NAS:} This method employs an intelligent agent (typically an RNN controller or Transformer) that learns a policy for sampling efficient architectures through a trial-and-feedback mechanism. In each upper level iteration, the RNN controller, also referred to as RNN-Architecture Sampler, generates a suitable architecture sequentially as per Markov Decision Process (MDP), wherein any architecture $\textsf{A}^i$ is represented as the list of actions $[A_1,A_2,\dots,A_T]$ (i.e., $\textsf{A}^i=[A_1,A_2,\dots,A_T]$). The lower level task involves implementing architecture $\textsf{A}^i$ and training the respective model using standard gradient-based method to minimize the training loss. This is a trial part of the trial-and-feedback mechanism. On the feedback side, the performance of trained model is evaluated on the validation set, and the performance result serves as a reward signal for the RNN-controller. Later, the reward signal information is used in policy gradient techniques to update the parameters of RNN controller such that the expected reward is maximized. Over the iterations, this procedure gradually refines the sampling distribution and directs the search toward more effective architectures. The recent surveys on the reinforcement learning based NAS methods (i.e., Automated Reinforcement Learning (AutoRL)) are available in \cite{jaafra2019reinforcement,parker2022automated}.
\end{enumerate}
The studies based on the methods discussed above are provided in Table \ref{tab:search_stra_NAS}.
\begin{table}[H] 
\centering
\caption{Search approaches followed in the sampling-based NAS methods} \label{tab:search_stra_NAS}
\begin{tabular}{|>{\raggedright\arraybackslash}p{5.8 cm}|
                >{\raggedright\arraybackslash}p{6 cm}|}
\hline
\textbf{Search Approach} & \textbf{Sample Studies} \\
\hline

Grid search & 
\cite{ogunsanya2023grid, dhilsath2021hyperparameter, radzi2021hyperparameter, ippolito2022hyperparameter, zahedi2021search, shekar2019grid}\\
\hline

Random search & 
\cite{dhilsath2021hyperparameter, ippolito2022hyperparameter, zahedi2021search, bergstra2012random, bergstra2011algorithms}\\
\hline

Bayesian optimization based search & 
\cite{radzi2021hyperparameter, ippolito2022hyperparameter, bergstra2011algorithms, eggensperger2013towards, snoek2012practical, klein2017fast}\\
\hline

Evolutionary computation based search & 
\cite{ippolito2022hyperparameter, alibrahim2021hyperparameter, lorenzo2017particle, byla2019deepswarm, suganuma2017genetic, chu2021fast, maziarz2018evolutionary, guo2020efficient, fan2022hybrid}\\
\hline

Reinforcement learning based search & 
\cite{zoph2016neural, chu2020multi, mohan2023autorl, hsu2018monas, gao2019graphnas}\\
\hline

\end{tabular}
\end{table}

\noindent\textbf{Bilevel Theory-based NAS Methods}: 
This class of methods formulates the search as a nested optimization problem. These methods are mostly based on the classical approach (Section \ref{sec:KKT} - Section \ref{sec:trust_region}) and they are built upon formal optimization theories, providing a systematic framework for NAS.

\begin{enumerate}[label=\arabic*)]
    \item \textbf{KKT Reformulation / MPEC-based NAS:} These methods reformulate the bilevel problem using KKT-conditions (\ref{eq:KKT_start} - \ref{eq:KKT_end}) or as a Mathematical Program with Equilibrium Constraints (MPEC) and use appropriate optimization methods to solve the reformulated problem.\vspace{2mm}

    \item \textbf{Hypergradient-based NAS:} These methods compute the hypergradients \cite{maclaurin2015gradient,sinha2025linear,shi2024doublemomentummethodlowerlevel}, i.e., derivatives of the validation loss with respect to architecture parameters, using the reverse mode or implicit differentiation. In some sense, it is similar to the KKT reformulation. DARTS (Differentiable Architecture Search) \cite{liu2018darts} and its variants fall into this category as it makes the architecture space continuous and differentiable.\vspace{2mm}

    \item \textbf{Penalty-based NAS:} These methods first reduce the bilevel problem to a single-level formulation and then relax the problem using penalty function or augmented Lagrangian function approach \cite{sinha2024gradient, shi2021improved}. The optimality of the lower level problem is ensured by adding regularization terms into the upper level objective.\vspace{2mm}

    \item \textbf{Surrogate Approximation-based NAS:} These methods approximate the $\Psi$-mapping or the $\varphi$-mapping using surrogate or meta-models. Approximating the $\Psi$-mapping emulates the model parameters in response to architecture change \cite{franceschi2018bilevel}. In \cite{sinha2024gradient}, $\varphi$-mapping is approximated using Kriging and then the single-level problem is solved using the Penalty-based approach, discussed above. \vspace{2mm}

    \item \textbf{Trust-Region-based NAS:} These techniques iteratively solve a series of approximate subproblems within a dynamically updated trust-region \cite{mackay2019self}. These methods ensure stable convergence in scenarios where small changes in hyperparameters lead to significant variations in model performance.\vspace{2mm}
\end{enumerate}    

\section{Conclusions and Future Research Directions}\label{sec:conclusions}
This paper presents a comprehensive review of bilevel optimization, covering fundamental principles and solution methods from both classical and evolutionary approaches. The nested structure of bilevel problems, where lower level optimization problem acts as a constraint of upper level optimization problem, often leads to non-convex and disconnected feasible region. These inherent complexities make bilevel optimization a challenging class of optimization problems. The methods from the classical approach are based on the rigorous mathematical optimization theories and they are suitable for bilevel problems following certain regularity conditions (i.e., problems consisting of mathematically well-behaving objective functions and constraints). The core mechanisms of such classical methods, including single-level reduction methods, descent methods, penalty function methods, and trust-region methods, are discussed. The evolutionary class of methods employs the nature-inspired or intelligence-based strategies to effectively explore the complex feasible region characterized by non-linearity, non-convexity, discreteness, discontinuity, etc. The use of evolutionary method along with classical method, a hybrid approach, is also discussed in the context of bilevel problem solving. Later, the mixed integer and multiobjective scenarios, which need further attention in bilevel optimization field, are also addressed. Two recent applications of bilevel approach, i.e., bilevel optimization-based decomposition and neural architecture search for machine learning algorithms, are discussed in detail. These new applications have the potential to impact the areas of optimization and machine learning significantly. Interestingly, the first application supports the automation of optimization problem solving, and the second application supports the automation of neural architecture design.

Although significant progress has been observed in recent years, bilevel optimization continues to be in a developmental stage, offering substantial scope for both computational and theoretical contributions.
Apart from that, discrete and multiobjective bilevel optimization topics have received limited attention in the bilevel research community. From a computational resource perspective, there is potential for using distributed computing platforms to effectively handle large-scale bilevel problems. 
Bilevel optimization remains an active research area, and the development of increasingly efficient algorithms is driving a shift toward more application-oriented studies. Simultaneously, emerging practical challenges are giving rise to newer types of practical bilevel problems. As a result, a broader spectrum of real-world applications is expected to surface in the near future.



\end{document}